\documentclass[a4paper]{article}
\usepackage{amssymb,amsmath}
\newtheorem{theo}{\sc Theorem}[section]
\newtheorem{prop}[theo]{\sc Proposition}
\newtheorem{rema}[theo]{\sc Remark}
\newcommand{\Z}{{\mathbb{Z}}}
\newcommand{\Q}{{\mathbb{Q}}}
\newcommand{\D}{{\mathbf{D}}}
\newcommand{\E}{{\mathbf{E}}}
\newcommand{\F}{{\mathbf{F}}}
\newcommand{\T}{{\mathbf{T}}}
\newcommand{\KE}{{\mathbf{e}}}
\newcommand{\KF}{{\mathbf{f}}}
\newcommand{\Sp}{{\mathrm{Sp}}}
\newcommand{\SL}{{\mathrm{SL}}}
\newcommand{\Ta}{{\mathbf{T}(p)}}
\newcommand{\Tb}{{\mathbf{T}_1(p^2)}}
\newcommand{\Tc}{{\mathbf{T}_2(p^2)}}
\newcommand{\Td}{{\mathbf{T}_3(p^2)}}
\newcommand{\Te}{{[\mathbf{p}]}}

\begin{document}

\title{Explicit Shimura's conjecture for $\Sp_4$}

\author{Kirill Vankov\\
\small Institut Fourier, Universit\'e Grenoble 1\\
\small UFR de Math\'ematiques, UMR 5582\\
\small BP 74, 38402 Saint-Martin d'H\`eres Cedex\\
\small France\\
\small e-mail: kvankov$@$fourier.ujf-grenoble.fr\\
\small phone: +33 476514656
}
\date{}
\maketitle

\begin{abstract}
Shimura's conjectire (1963) concerns the rationality of the generating  series for Hecke operators for the symplectic group of genus $g$. This conjecture  wes proved by Andrianov for arbitrary genus $g$. For genus $g=4$, we explicify the rational fraction in this conjecture.  Using formulas for images of double cosets, we first compute the sum of the generating series under the Satake spherical map, which is a rational fraction with polynomial coefficients.  Then we recover the coefficients of this fraction as elements of the Hecke algebra using polynomial representation of basic Hecke operators under spherical map.  Numerical examples of these fractions for special choice of Satake parameters are given.
\end{abstract}

\tableofcontents

\section{Introduction}
\label{sec:intro}

Let $p$ be a prime. We consider the symplectic group $\Sp_n$ of genus $n$, and let
\begin{equation*}
\{\Ta,\Tb,\dotsc,\T_{n-1}(p^2),\Te_n\}
\end{equation*}
be $n+1$ generators of the Hecke ring over $\Z$ for $\Sp_n\subset\mathrm{GL}_{2n}$. Let $D_p(X)$ denote the generating power series of Hecke operators
\begin{equation}\label{eq:Dp(X)}
\D_{p}(X)=\sum_{\delta=0}^{\infty} \T(p^\delta )X^\delta \text{ .}
\end{equation}

\

The result presented in this article provides a complement to the solution of Shimura's conjecture of rationality of generating Hecke power series stated in \cite{Sh63} at p.~825 as follows:
\begin{quote}
``In general, it is plausible that $\D_p(X)=\E(X)/\F(X)$ with polynomials $\E(X)$  
and $\F(X)$ in $X$ with integral coefficients of degree $2^n-2$ and $2^n$, respectively'' 
\end{quote}
(\emph{i.e. with coefficients in Hecke algebra} $\mathcal{L}_\Z=\Z[\Ta,\Tb,\dotsc,\T_{n-1}(p^2),\Te_n]$).

The existence of a rational representation $\E(X)/\F(X)$ was proved by Andrianov in \cite{An67,An68,An69}
for arbitrary genus $n$.  For genus 1 and 2 the results were given by Hecke and Shimura (\cite{He59}, \cite{Sh71}, Theorem 3.21, and \cite{Sh63}, Theorem 2):
\begin{align*}
\D&_{p}^{(1)}(X)=\dfrac{1}{1-\Ta X+p\Te_1X^2} \;,\\
&\\
\D&_{p}^{(2)}(X)=\\
&\dfrac{1-p^2\Te_2X^2}{1-\Ta X+p\,(\Tb+(p^2+1)\Te_2)X^2-p^3\Te_2\Ta X^3+p^6\Te_2^2X^4}\;.\\
&
\end{align*}
Andrianov obtained the expression for genus 3 using the multiplication table of Hecke operators in \cite{An67}.  No explicit results for higher genus were known due to the enormous complexity of Hecke algebra manipulations.  Recently the author together with Panchishkin developed a formal calculus approach using a computer.  We were able to compute more directly the generating series in Shimura's conjecture for genus 3 (see \cite{PaVaSp3}), and then to explore the case of genus 4.  Here is the result for genus 3, where coefficients in $p$ are factorized into irreducible polynomials:
\begin{align}
\label{eq:D(n=3)}
&\D_{p}^{(3)}(X)=\dfrac{\E_3(X)}{\F_3(X)} \text{ , where } \E_3(X), \F_3(X) \in \mathcal{L}_{\Z[X]} \text{ and}
\\\nonumber
&\E_3(X)=\\\nonumber
& \quad 1-p^2\left(\Tc+(p^2-p+1)(p^2+p+1)\Te_3\right)X^2 + p^4(p+1)\Te_3\Ta X^3 \\\nonumber
& \quad -p^7\Te_3\left(\Tc+(p^2-p+1)(p^2+p+1)\Te_3\right) X^4 + p^{15}\Te_3^3\,X^6 \;,
\end{align}
\begin{align}\nonumber
&\F_3(X) = 1 - \Ta X \\\nonumber
& \quad +p\left(\Tb+(p^2+1)\Tc+(p^2+1)^2\Te_3\right) X^2 \\\nonumber
& \quad -p^3\left(\Tc+\Te_3\right)\Ta X^3 \\\nonumber
& \quad +p^6\big(\Tc+\Te_3(\Ta^2-2p\Tb-2(p-1)\Tc\\\nonumber
&\qquad-(p^2+2p-1)(p^2-p+1)(p^2+p+1)\Te_3)\big)X^4 \\\nonumber
& \quad -p^9\Te_3\left(\Tc+\Te_3\right)\Ta X^5 \\\nonumber
& \quad +p^{13}\Te_3^2\left(\Tb+(p^2+1)\Tc+(p^2+1)^2\Te_3\right) X^6\,.
\end{align}

In this article we describe the application of formal calculus for genus 4 (case $\Sp_4$).  We present both numerator and denominator polynomials expressed in terms of Hecke operators in section \ref{sec:Sp4}.  Notations of the article and some useful facts of Hecke algebras are given in section \ref{sec:hecke}.  In section \ref{sec:map} we define the Satake mapping of Hecke algebra to the symmetrical polynomial ring, which is referred by Andrianov and Zhuravlev in \cite{AnZh95} as the spherical map.  Then we discuss in details the method of obtaining the main result in section \ref{sec:main}.  Finally we give some interesting properties of obtained polynomials in section \ref{sec:rema}.

Generating series of a type
\begin{equation*}
\sum_{m=1}^{\infty} \lambda_f(m)m^{-s}=\prod_{p\,\mathrm{primes}}\sum_{\delta=0}^{\infty} \lambda_f(p^\delta )p^{-\delta s}
\end{equation*}
are used as a classical method to produce $L$-functions for an algebraic group $G$ over $\Q$, where $\lambda_f(m)$ are the eigenvalues of Hecke operators on an automorphic form $f$ on $G$.  Hence these series and related congruences are of number-theoretic interest.  Particularly, we study here the generating series of Hecke operators $\T(m)$ for the symplectic group $\Sp_4$ and $\lambda_f(m)=\lambda_f(\T(m))$.

The explicit knowledge of the sum of the generating series of Hecke operators
$$
\D_p(X)=\sum_{\delta=0}^{\infty} \T(p^\delta )X^\delta=\E(X)/\F(X)
$$
gives a relation  between the Hecke eigenvalues and the Fourier coefficients of a Hecke eigenform $f$.  This link is needed for constructing an analytic continuation of $L$-function on $\Sp_n$, which was done by Andrianov for $\Sp_2$ in \cite{An74}.  An approach for constructing an analytic continuation of the spinor $L$-function on $\Sp_3$ was indicated by Panchishkin at the talk on seminar Groupes R\'eductifs et Formes Automorphes in the Institut de Math\'ematiques de Jussieu (\cite{PaGRFA}).

Similar technique of a symbolic computation can be used to discover other interesting identities between Hecke operators, between their eigenvalues, relations to Fourier coefficients of modular forms of higher degree.  In \cite{PaVaRnk} the author together with Panchishkin study the analogue of Rankin's Lemma of higher genus.

\section{The Explicit Formula for $\Sp_4$}
\label{sec:Sp4}

\begin{theo}
\label{th:main}
For genus $g=4$ the summation of Hecke power series $\D_{p}(X)$ resolves explicitly to the following rational polynomial presentation:
$$
\D_{p}^{(4)}(X)=\sum_{\delta=0}^{\infty} \T(p^\delta )X^\delta=\dfrac{\E_4(X)}{\F_4(X)}\;,
$$
where $\E_4(X)=\sum_{k=0}^{14}\KE_kX^k$ is the polynomial of degree $14$ and $\F_4(X)=\sum_{k=0}^{16}\KF_kX^k$ is the polynomial of degree $16$ with the coefficients $\KE_k$ and $\KF_k$ listed below:

\begin{align*}
\KE_{0} &= 1
\,,\\
\KE_{1} &= 0
\,,\\
\KE_{2} &= -p^2 (
 (p^8+p^6+2p^4+2p^2+1) \Te
+(p^2+p+1) (p^2-p+1)   \Td
+                      \Tc
)
\,,\\
\KE_{3} &= p^4 (p+1) (
 (p^2+1) (p^3-p^2+1) \Te
+                    \Td
)\Ta
\,,\\
\KE_{4} &= p^7 (
 (p^2+p+1) (p^2-p+1) (p^8+3p^7+p^5+2p^3+p-1) \Te^2     \\ & \quad
+(p^2+p+1) (p^2-p+1) (2p^3+p-2)              \Td   \Te
-(p^2+p+1) (p^2-p+1)                         \Td^2     \\ & \quad
+(2p^5+2p^3+p^2+p-1)                         \Tc   \Te
-                                            \Tc   \Td
+p (p^2+p+1)                                 \Tb   \Te \\ & \quad
-(p^2+p+1)                                   \Ta^2 \Te
)
\,,\\
\KE_{5} &= - p^{10} (p+1) (
 (p^2+1) (p^7-p^6-p^2-1)   \Te
-(p^2+1)                   \Td
-                          \Tc
) \Ta \Te
\,,\\
\KE_{6} &= p^{14} (
 (p^{16}-p^{15}-2p^{14}-3p^{12}-5p^{10}-8p^8+p^7-8p^6-5p^4-4p^2-1) \Te^3           \\ & \quad
+(p^{12}-p^{10}-p^9-5p^8+2p^7-7p^6-6p^4-8p^2+p-2)                  \Td   \Te^2     \\ & \quad
+(p^7-p^4-4p^2+2p-1)                                               \Td^2 \Te
+p                                                                 \Td^3           \\ & \quad
-(2p^8+3p^6+p^4-p^3+3p^2+p+1)                                      \Tc   \Te^2     \\ & \quad
-(3p^2+p+1)                                                        \Tc   \Td   \Te
-(p^6-p^3+1)                                                       \Tb   \Te^2
-                                                                  \Tb   \Td   \Te \\ & \quad
+p^2 (p^3+p-1)                                                     \Ta^2 \Te^2
)
\,,\\
\KE_{7} &= - p^{19} (p-1) (p+1) (
 (p^2+p+1) (p^2-p+1) (p^2+1)   \Te  \\ & \quad
+(p^2+p+1) (p^2-p+1)           \Td
+                              \Tc
) \Ta \Te^2
\,,\\
\KE_{8} &= - p^{24} (
 (p^{16}-3 p^{12}-3 p^{10}-p^9-9 p^8-8 p^6-7 p^4-5 p^2+p-1)     \Te^3       \\ & \quad
+(p^{10}-p^9-4 p^8-6 p^6-8 p^4-9 p^2+3 p-2)             \Td     \Te^2       \\ & \quad
-(p^4+4 p^2-3 p+1)                                      \Td^2   \Te
+p                                                      \Td^3               \\ & \quad
-(p^8+2 p^6-p^5+2 p^4+p^3+4 p^2+1)                      \Tc     \Te^2       \\ & \quad
-(p^3+3 p^2+1)                                          \Tc     \Td   \Te
+(p^5-p^2-1)                                            \Tb     \Te^2       \\ & \quad
-                                                       \Tb     \Td   \Te
-p (p^3-p^2-1)                                          \Ta^2   \Te^2
) \Te
\,,\\
\KE_{9} &= p^{29} (p+1) (
 (p^2+1) (p^5-2 p^4-1)   \Te
-(p^4+1)                 \Td
-                        \Tc
) \Ta \Te^3
\,,
\\
\KE_{10} &= - p^{35} (
 (p^2-p+1) (p^2+p+1) (p^8+2 p^7+p^5+3 p^3+p-1)   \Te^2         \\ & \quad
-(p^2-p+1) (p^2+p+1) (p^5-3 p^3-p+2)             \Td     \Te   \\ & \quad
-(p^2-p+1) (p^2+p+1)                             \Td^2
+(p^5+3 p^3+p^2+p-1)                             \Tc     \Te   \\ & \quad
-                                                \Tc     \Td
+p (p^2+p+1)                                     \Tb     \Te
-(p^2+p+1)                                       \Ta^2   \Te
) \Te^3
\,,
\end{align*}
\begin{align*}
\KE_{11} &= - p^{41} (p+1) (
(p^2+1) (p^3-p^2+1)   \Te
+                     \Td
) \Ta \Te^4
\,,\\
\KE_{12} &= p^{48} (
 (2 p^6+2 p^4+2 p^2+1)   \Te
+(p^2-p+1) (p^2+p+1)     \Td
+                        \Tc
) \Te^5
\,,\\
\KE_{13} &= 0
\,,\\
\KE_{14} &= - p^{64} \Te^7\,,
\end{align*}
\begin{align*}
\KF_{0} &= 1\,,\quad\KF_{1} = -\Ta
\,,\\
\KF_{2} &= p (
 (p^8+2 p^6+2 p^4+2 p^2+1)  \Te
+(p^4+2 p^2+1)              \Td
+(p^2+1)                    \Tc
+                           \Tb
)
\,,\\ 
\KF_{3} &= p^3 (
 (p^7-p^6-p^4-p^2-1)  \Te
-                     \Td
-                     \Tc
) \Ta
\,,\\
\KF_{4} &= -p^6 (
(p^4+1)^2 (p^6-p^4+2 p^3-2 p^2+2 p-1)    \Te^2        \\ & \quad
+2 (p^2-p+1) (p^4+1) (p^4+2 p^3+p^2+p-1) \Td \Te      \\ & \quad
+(p^2-p+1) (p^2+p+1) (p^2+2 p-1)         \Td^2
-2 (p^2-p+1) (p^4+1)                     \Tc \Te      \\ & \quad
+2 (p-1)                                 \Tc \Td
-                                        \Tc^2
+2 p (p^4+1)                             \Tb \Te
+2 p                                     \Tb \Td      \\ & \quad
-                                        \Ta^2 \Te
-                                        \Ta^2 \Td
)
\,,
\\
\KF_{5} &= p^9 (
 (p^{11}+p^{10}+4 p^8+2 p^7+3 p^6+3 p^4+2 p^2-1)  \Te^2     \\ & \quad
+(2 p^7+2 p^6+3 p^4+2 p^2-2)                  \Td \Te
-                                             \Td^2
+(p^4+3 p^2-1)                                \Tc \Te       \\ & \quad
-                                             \Tc \Td
+3 p^2                                        \Tb \Te
-p                                            \Ta^2 \Te
) \Ta
\,,\\
\KF_{6} &= -p^{13} (
 (p^4+1) (p^2+1)^2 (2 p^8+2 p^7+2 p^5+2 p^3+2 p-1) \Te^3          \\ & \quad
+(p^2+1)^2 (2 p^8+2 p^7+4 p^5-2 p^4+2 p^3+4 p-3)   \Td \Te^2      \\ & \quad
-(p^2+1)^2 (p^4-2 p+3)                             \Td^2 \Te
-(p^2+1)^2                                         \Td^3          \\ & \quad
+(p^2+1) (2 p^8+4 p^7+4 p^5+4 p^3+4 p-1)           \Tc \Te^2      \\ & \quad
+2 (p^2+1) (p^3+2 p-1)                             \Tc \Td \Te
-(p^2+1)                                           \Tc \Td^2      \\ & \quad
+2 p (p^2+1)                                       \Tc^2 \Te
+(2 p^8+2 p^7+2 p^5+2 p^3+2 p-1)                   \Tb \Te^2      \\ & \quad
+2 (p-1)                                           \Tb \Td \Te
-                                                  \Tb \Td^2 
+2 p                                               \Tb \Tc \Te    \\ & \quad
-(p^2+1) (p^5+p^4-p^3+1)                           \Ta^2 \Te^2
+(p-1) (p^2+p+1)                                   \Ta^2 \Td \Te
-                                                  \Ta^2 \Tc \Te
)
\,,
\\
\KF_{7} &= -p^{17} (
 (2 p^{13}+p^{12}+3 p^{10}+p^9+p^8+2 p^6-p^5+p^4-p^2+p+1) \Te^3   \\ & \quad
+(p^9+2 p^6-2 p^5+2 p^4-2 p^2+3 p+2)                \Td \Te^2     \\ & \quad
-(p^5-p^4+p^2-3 p-1)                                \Td^2 \Te
+p                                                  \Td^3
+(p^6+2 p^4-2 p^2+1)                                \Tc \Te^2     \\ & \quad
-(2 p^2-1)                                          \Tc \Td \Te
+(2 p^4+1)                                          \Tb \Te^2
+                                                   \Tb \Td \Te
-p^3                                                \Ta^2 \Te^2
) \Ta
\,,
\end{align*}
\begin{align*}
\KF_{8} &= p^{22} (
 (p^{18}+4 p^{17}+3 p^{16}+8 p^{15}+12 p^{14}+8 p^{13}+14 p^{12}+12 p^{11}+20 p^{10}          \\ & \qquad
                                   +4 p^9+20 p^8+16 p^6+10 p^4-4 p^3+5 p^2+1) \Te^4           \\ & \quad
+2 (2 p^{10}+2 p^9+p^8+6 p^7+4 p^6+8 p^4-6 p^3+6 p^2+1) (p^4+1)                 \Td \Te^3     \\ & \quad
+(p^8+4 p^6+8 p^4-12 p^3+10 p^2+1)                                            \Td^2 \Te^2
-4 p^2 (p-1)                                                                  \Td^3 \Te
+p^2                                                                          \Td^4           \\ & \quad
+2 (2 p^7+3 p^6+2 p^5+5 p^4-2 p^3+3 p^2+1) (p^4+1)                            \Tc \Te^3       \\ & \quad
+(2 p^6+4 p^5+10 p^4-8 p^3+6 p^2+2)                                           \Tc \Td \Te^2
-4 p^3                                                                        \Tc \Td^2 \Te   \\ & \quad
+(3 p^4+2 p^2+1)                                                              \Tc^2 \Te^2
+2 (2 p^5+p^4+2 p^2+1) (p^4+1)                                                \Tb \Te^3       \\ & \quad
+(4 p^5+2 p^4+4 p^2+2)                                                        \Tb \Td \Te^2
+2 (p^2+1)                                                                    \Tb \Tc \Te^2   \\ & \quad
+                                                                             \Tb^2 \Te^2
-(p^8+2 p^7+2 p^5+2 p^4+2 p^3+2 p-1)                                          \Ta^2 \Te^3     \\ & \quad
-2 (p^4+p-1)                                                                  \Ta^2 \Td \Te^2
+                                                                             \Ta^2 \Td^2 \Te
-2 p                                                                          \Ta^2 \Tc \Te^2
)
\end{align*}
and the higher degree coefficients are the same as the lower ones multiplied by $p^{10}\Te$ in the correspondent degree:
\begin{align*}
\KF_{9}  &= \KF_{7} \cdot p^{10} \Te   ,&
\KF_{10} &= \KF_{6} \cdot p^{20} \Te^2 ,&
\KF_{11} &= \KF_{5} \cdot p^{30} \Te^3 ,&
\KF_{12} &= \KF_{4} \cdot p^{40} \Te^4 ,&
\\
\KF_{13} &= \KF_{3} \cdot p^{50} \Te^5 ,&
\KF_{14} &= \KF_{2} \cdot p^{60} \Te^6 ,&
\KF_{15} &= \KF_{1} \cdot p^{70} \Te^7 ,&
\KF_{16} &= \KF_{0} \cdot p^{80} \Te^8 .&
\end{align*}
\end{theo}


The proof of this result is based on application of spherical map in order to carry all calculations in the ordinary polynomial ring instead of Hecke algebra.  Using formal calculus on a computer it is possible to find the explicit symmetrical polynomial solution for the image $\Omega(\D_p(X))$.  Similarly we can find the images of basis Hecke operators, which we use to compose and resolve a linear system of undetermined coefficients and discover the desired expression in terms of Hecke operators.

\section{The Hecke Algebras}
\label{sec:hecke}

In this section we describe the notations used in the article.  Most definitions are taken from \cite{AnZh95}, where the detailed theory of Hecke rings is given in Chapter 3 (see also \cite{An87}).

\subsection{Hecke algebra for $\Sp_n$}

Consider the group of positive symplectic similitudes
\begin{align*}
&\mathrm{S}=\mathrm{S}^n=\mathrm{GSp}_n^+(\Q)=\{M\in \mathrm{M}_{2n}(\Q)\;\vert\;{}^tMJ_nM=\mu(M)J_n, \mu(M)>0\}\;,\\
&\text{where } J_n = \left( \begin{array}{cc}\mathbf{0}_n& \mathbf{I}_n \\ -\mathbf{I}_n& \mathbf{0}_n\end{array}\right).
\end{align*}

For the Siegel modular group $\Gamma=\Sp_n(\Z)\subset \SL_{2n}(\Z)$ of genus $n$ consider the double cosets
$$
\big(M\big)=\Gamma M\Gamma\subset\mathrm{S}\,,
$$
and the Hecke operators
$$
\T(\mu) = \sum_{M\in \mathrm{SD}_n(\mu)}\big(M\big)\,,
$$
where $M$ runs through the following integer matrices
$$
\mathrm{SD}_n(\mu)=\{\mathrm{diag}(d_1,\dotsc,d_n;e_1,\dotsc,e_n)\}\,,
$$
where $d_1|\dotsb|d_n|e_n|\dotsb|e_1,\;d_i,e_j>0,\;d_ie_i=\mu=\mu(M)$.  Let us use the notation for the Hecke operators
$$
\T(d_1,\dotsc,d_n;e_1,\dotsc,e_n)=\big(\mathrm{diag}(d_1,\dotsc,d_n;e_1,\dotsc,e_n)\big)\,.
$$
In particular we have the following $n+1$ basis Hecke operators
\begin{equation}
\label{eq:T-basis}
\begin{split}
&\Ta=\T(\underbrace{1,\dotsc,1}_{n},\underbrace{p,\dots,p}_{n})\,,\\
&\T_i(p^2)=\T(\underbrace{1,\dotsc,1}_{n-i},\underbrace{p,\dotsc,p}_{i},\underbrace{p^2,\dotsc,p^2}_{n-i},\underbrace{p,\dotsc,p}_{i})\,,\,i=1,\dotsc,n\,,
\end{split}
\end{equation}
generating the Hecke algebra over $\Z$:
$$
\mathcal{L}_{n,\Z}=\Z[\Ta,\Tb,\dotsc,\T_n(p^2)]\,.
$$
We denote as $\Te_n$ (or just $\Te$ if the context of $n$ is declared) the scalar matrix Hecke operator $\Te=\T_n(p^2)=\T(\underbrace{p,\dotsc,p}_{2n})=p\mathbf{I}_{2n}$ .

\subsection{Operation of multiplication in Hecke algebras}

In order to define an operation of multiplication (in abstract Hecke algebra) we consider without lost of generality for any subgroup $\Gamma$ of semigroup $S$ a vector space over $\Q$ generated by all left cosets $\Gamma M$ (just as a formal base)
$$
L_{\Q}(\Gamma,S)=\Big\lbrace\sum_ja_j(\Gamma M_j)\;\vert\;a_j\in\Q\Big\rbrace\;.
$$
Only for double cosets $\big(M\big)=\Gamma M\Gamma\subset\mathrm{S}$ that can be presented as a finite union of disjoint left cosets 
$$
\big(M\big)=\bigcup_{j=1}^K\Gamma M_j\;,
$$
we denote
$$
\big(M\big)=\sum_{M_j\in\Gamma\setminus\Gamma M\Gamma}(\Gamma M_j)
$$
and consider an abstract Hecke algebra $\mathcal{L}_{\Q}(\Gamma,S)=L_{\Q}(\Gamma,S)^{\Gamma}$ as a vector space $L_{\Q}(\Gamma,S)^{\Gamma}$ for fixed $\Gamma$.  Any nonzero element $t\in\mathcal{L}$ can be written in the form $t=\sum_{j=1}^Ka_j\big(\Gamma M_j\big)$.  Hence, the multiplication is well defined as
$$
\Big(\sum_ja_j\big(\Gamma M_j\big)\Big)\Big(\sum_{j'}a^\prime_{j'}\big(\Gamma M_{j'}\big)\Big)=\sum_{j,j'}a_ja^\prime_{j'}\big(\Gamma M_jM_{j'}\big),\;a_j,a^\prime_{j'}\in\Q\;.
$$

\subsection{Hecke algebra for $\mathrm{GL}_n$}

Further, in order to define a mapping to the polynomial ring, we need to introduce a Hecke algebra for general linear group.  Let $G=\mathrm{GL}_n(\Q)$ and $\Lambda=\mathrm{GL}_n(\Z)$.  We note corresponding Hecke algebra as $\mathcal{H}_{\Q}(\Lambda,G)=L_{\Q}(\Lambda,G)^{\Lambda}$.  This algebra is generated by only $n$ basis operators
\begin{equation}
\label{eq:pi}
\pi_i(p)=\big(\mathrm{diag}(\underbrace{1,\dotsc,1}_{n-i},\underbrace{p,\dotsc,p}_{i})\big)\,,\;1\leqslant i\leqslant n\;.
\end{equation}
Recall that every left coset $\Lambda g$ ($g\in G$) has a representative of the form
\begin{equation}
\label{eq:g-rep}
\begin{pmatrix}
  p^{\delta_1} & c_{12}       & \dotsb & c_{1n}     \\
  0            & p^{\delta_2} & \dotsb & c_{2n}      \\
  \dotsb       & \dotsb       & \dotsb & \dotsb      \\
  0            & 0            & \dotsb & p^{\delta_n}
\end{pmatrix}
\,,\text{ where }\delta_1,\dotsc,\delta_n\in\Z\,.
\end{equation}
The arbitrary element $t\in\mathcal{H}_{\Q}(\Lambda,G)$ is composed as a finite linear combination of left cosets $\Lambda g_j$:
\begin{equation}
\label{eq:t}
t=\sum_ja_j\big(\Lambda g_j\big)\,.
\end{equation}

We need to introduce another element related to the product of generators $\pi_i(p)$, which will be later used for spherical map definition.  Let $\pi_{\alpha\beta}$ be a double coset defined by
\begin{equation}
\label{eq:pi_ab}
\pi_{\alpha\beta} = \pi_{\alpha\beta}^n(p) = 
\left(\begin{pmatrix}
  \mathbf{I}_{n-\alpha-\beta} & 0 & 0 \\
  0 & p\,\mathbf{I}_\alpha & 0 \\
  0 & 0 & p^2\mathbf{I}_\beta
\end{pmatrix}\right)\;.
\end{equation}
The double coset expansion of the product in the Hecke algebra $\mathcal{H}_\Q$ of two generators $\pi_i$ and $\pi_j$, where $1\leqslant i,j\leqslant n$, has the form
$$
\pi_i\pi_j=\sum_{\substack{0\leqslant a\leqslant n-j\\0\leqslant b\leqslant j\\a+b=i}}\dfrac{\varphi_{a+j-b}}{\varphi_a\varphi_{j-b}}\pi_{a+j-b,b}\;,
$$
where
\begin{align}
\label{eq:varphi}
&\varphi_i=\varphi_i(x)=(x-1)(x^2-1)\dotsb(x^i-1)\text{ for }i\geqslant 1\\
&\nonumber\text{and }\varphi_0(x)=1\;.
\end{align}

\section{The Spherical Map}
\label{sec:map}

There are several methods to construct a mapping from a Hecke algebra to a polynomial ring.  We use as a base the book of Andrianov and Zhuravlev \cite{AnZh95}, where the description for general linear and symplectic groups is given in terms of right cosets of the double cosets, which generate the Hecke algebra.  This isomorphism plays a key role in our calculations.  It allows to carry all computation in the polynomial ring where a multiplication is much more straightforward than the product of double cosets.

\subsection{The spherical map in general linear group case}

The spherical map for the Hecke algebra of general linear group is defined for fixed representative of a left coset of the form (\ref{eq:g-rep}) as
\begin{equation*}
\omega(\big(\Lambda g\big)) = \prod_{i=1}^n(x_ip^{-i})^{\delta_i}
\end{equation*}
and for an arbitrary element (\ref{eq:t}) we have
\begin{equation*}
\omega(t) = \sum_ja_j\omega(\big(\Lambda g_j\big))\,.
\end{equation*}
This definition is unique due to the fact that the diagonal $(p^{\delta_1},\dotsc,p^{\delta_n})$ in (\ref{eq:g-rep}) is uniquely determined by the left coset.

Lemma 2.21 of chapter 3 in \cite{AnZh95} gives the images of the basis elements (\ref{eq:pi}) of Hecke algebra for general linear group:
\begin{equation*}
\omega(\pi_i(p))=p^{-\langle i\rangle}s_i(x_1,\dotsc,x_n)\quad(1\leqslant i\leqslant n),
\end{equation*}
where
\begin{equation*}
s_i(x_1,\dotsc,x_n)=\sum_{1\leqslant\alpha_1<\dotsb<\alpha_i\leqslant n}x_{\alpha_1}\dotsb x_{\alpha_i}
\end{equation*}
is the $i$--th elementary symmetric polynomial.

\subsection{The spherical map in symplectic group case}

The definition of the spherical map in case of symplectic group is more complicated and based on the case of general linear group.  Consider an arbitrary element (double coset) $T\in\mathcal{L}_{n,\Z}$ as a finite linear combination of left cosets:
\begin{equation*}
T=\sum_jb_j\big(\Gamma M_j\big)\,,\text{ with }\mu(M_j)=p^{\delta_j}\,.
\end{equation*}
We choose the representative of a class in the form
\begin{equation*}
M_j=\begin{pmatrix}p^{\delta_j}D_j^*&*\\0&D_j\end{pmatrix}\,,
\end{equation*}
where $D_j^*={}^tD_j^{-1}$ and matrix $D_j$ is a triangular of the form
\begin{equation*}
D_j=\begin{pmatrix}
p^{\gamma_{1j}} & *               & \dotsb & \dotsb         \\
0               & p^{\gamma_{2j}} & *      & \dotsb         \\
\dotsb          & \dotsb          & \dotsb & \dotsb         \\
0               & \dotsb          & 0      & p^{\gamma_{nj}}
\end{pmatrix}\,.
\end{equation*}
We define the mapping as
\begin{equation*}
\Omega(T)=\sum_jb_jx_0^{\delta_j}\omega(\Lambda D_j)\,.
\end{equation*}

In particular on page 146 of \cite{AnZh95} we have the following formulas for basis Hecke operators:
\begin{align}
\label{eq:Ta}
&\Omega(\Ta)=\sum_{a=0}^nx_0s_a(x_1,\dotsc,x_n)=x_0\prod_{i=1}^n(1+x_i)\,,\\
\label{eq:T}
&\Omega(\T_i(p^2))=\sum_{\substack{a+b\leqslant n\\a\geqslant i}}
  p^{b(a+b+1)}l_p(a-i,a)x_0^2\omega(\pi_{a,b}(p))\,,
\end{align}
where the coefficient $l_p(r,a)$ is the number of $a\times a$ symmetric matrices of rank $r$ over the field of $p$ elements.  This coefficient is explicitly given by the recurrent formula (6.79) of \cite{AnZh95} on page 214 (Chapter 3, \S6):
\begin{equation*}
l_p(r,a)=l_p(r,r)\frac{\varphi_a(p)}{\varphi_r(p)\varphi_{a-r}(p)}\,,
\end{equation*}
where the function $\varphi_i(x)$ was already defined by (\ref{eq:varphi}).

Now we apply the above formulas to the power series (\ref{eq:Dp(X)}).  Following \cite{AnZh95} on page 150 we obtain the expression for the image of $\D_p(X)$:
\begin{align}
\label{eq:Omega(D)}
\Omega(&\D_p(X))=
\sum_{\delta=0}^{\infty}\Omega(\T(p^\delta))X^\delta =
\\\nonumber
&\sum_{\delta=0}^{\infty}\quad
 \sum_{0\leqslant\delta_1\leqslant\dotsb\leqslant\delta_n\leqslant\delta}\,
 p^{n\delta_1+(n-1)\delta_2+\dotsb+\delta_n}\,
 \omega(t(p^{\delta_1},\dotsc,p^{\delta_n}))\,(x_0X)^\delta\,,
\end{align}
where 
\begin{equation}
\label{eq:t(p)}
t(p^{\delta_1},\dotsc,p^{\delta_n})=\big(\mathrm{diag}(p^{\delta_1},\dotsc,p^{\delta_n})\big)\in\mathcal{H}_{\Q}
\end{equation}
is an element of the Hecke algebra for general linear group.

\subsection{Practical computation}

The algorithm was programmed and the results were computed using Maple system.  We found more practical and suitable for direct programming the formulas for spherical mapping in the article \cite{An70}.  Note that the notation $\Omega$ denotes in that article the spherical mapping of the Hecke algebra for general linear group.  It corresponds to our mapping $\omega$ defined above with substitution of all $x_i$ by $x_i/p$ for $i=1,\dotsc,n$.  Therefore we used the formula (1.7) on page 432 of \cite{An70} and then performed the substitution.  This formula gives direct expression for images of the elements $t$ of a type (\ref{eq:t(p)}) including images of $\pi_{\alpha\beta}(p)$ (\ref{eq:pi_ab}).  In our notation it can be written as
\begin{equation}
\label{eq:omega(t)}
\omega(t(p^{(\delta)})=p^{-\sum_i(n-i)\delta_i}\frac{Q(x)}{P^{(k)}(\frac{1}{p})}\,,
\end{equation}
where
\begin{align*}
&P^{(k)}(x)=\frac{\varphi_{k_1}(x)\dotsc\varphi_{k_t}(x)}{\varphi_1(x)^n}\,,\\
&Q(x)=\sum_{w\in \mathrm{S}_n}(wx)^{(\delta)}c(wx)\,,\\
&c(x)=\prod_{\alpha\in\Sigma}\frac{1-\frac{1}{p}(x)^{(\alpha)}}{1-(x)^{(\alpha)}}\,,
\end{align*}
function $\varphi(x)$ was defined by (\ref{eq:varphi}), the notation $(x)$ is used for $(x_1,x_2,\dotsc,x_n)$, then $(x)^{(\alpha)}\equiv x_1^{\alpha_1}x_2^{\alpha_2}\dotsb x_n^{\alpha_n}$, $(wx)^{(\delta)}\equiv x_{w(1)}^{\delta_1}x_{w(2)}^{\delta_2}\dotsb x_{w(n)}^{\delta_n}$. The set $\Sigma=\{(\alpha)\}=\{(\alpha_1,\alpha_2,\dotsc,\alpha_n)\}=\{\alpha_{ij}, 1\leqslant i<j\leqslant n\}$, where $\alpha_{ij}$ is defined by placing of $1$ and $-1$ within the set of $n$ zeros $\alpha_{ij}\equiv(\dotsc,0,1_{i},0,\dotsc,0,-1_{j},0,\dotsc)\in\Z_n$.  The element of Hecke algebra for the general linear group noted as $t(p^{(\delta)})$ is $t(p^{\delta_1},\dotsc,p^{\delta_n})$.  Numbers $(k)\equiv(k_1,\dotsc,k_t)$ denote the quantities of $t$ distinct elements in the set of integers $(\delta)=(\delta_1,\dotsc,\delta_n)$, that is the number $\delta_1$ occurs in $(\delta)$ exactly $k_1$ times, the next number following $\delta_1$ in the ordering of $(\delta)$ and distinct from $\delta_1$ appears there $k_2$ times, etc.  Note, that all $k_i>0$ and $k_1+\dotsb+k_t=n$.

For our computation we consider $n=4$. The set $\Sigma$ consists of 6 elements $\Sigma=\{(1,-1,0,0), (1,0,-1,0), (1,0,0,-1), (0,1,-1,0), (0,1,0,-1), (0,0,1,-1)\}$.  In consequence
\begin{equation*}
c(x)=\frac{(px_2-x_1)(px_3-x_1)(px_4-x_1)(px_3-x_2)(px_4-x_2)(px_4-x_3)}{p^6(x_2-x_1)(x_3-x_1)(x_4-x_1)(x_3-x_2)(x_4-x_2)(x_4-x_3)}\,.
\end{equation*}
In order to find images of Hecke operators (\ref{eq:T}) we need to apply the above formula (\ref{eq:omega(t)}) for all $t(p^{(\delta)})$ of a type (\ref{eq:pi_ab}), that is for the following $\{(\delta)\}=\{(0,0,0,1)$,  $(0,0,1,2)$, $(0,1,2,2)$, $(1,2,2,2)$, $(0,0,1,1)$, $(0,1,1,2)$, $(1,1,2,2)$, $(0,1,1,1)$, $(1,1,1,2)$, $(1,1,1,1)\}$.  Further we will see that in order to compute the series (\ref{eq:Dp(X)}) we need all $\{(\delta)\}$ with components up to 14.  We can dramatically reduce this set by taking the common degree of $p$ from $t(p^{\delta_1},\dotsc,p^{\delta_n})$ outside of the double coset matrix (multiplying the element by the corresponding degree of $p$).  Therefore, we need to compute just 680 primitive elements of the form $t(1,p^{\delta_2},p^{\delta_3},p^{\delta_4})$, where $0\leqslant\delta_2\leqslant\delta_3\leqslant\delta_4\leqslant 14$.

Here are some examples of values of these images:
\begin{align*}
&\omega(t(1,1,1,1))=1\,,\\
&\omega(t(1,1,1,p))=p^{-1}(x_1+x_2+x_3+x_4)\,,\\
&\omega(t(1,1,p,p))=p^{-3}(x_1x_2+x_1x_3+x_1x_4+x_2x_3+x_2x_4+x_3x_4)\,,\\
&\omega(t(1,p,p,p))=p^{-3}(x_1x_2x_3+x_1x_2x_4+x_1x_3x_4+x_2x_3x_4)\,,\\
&\omega(t(1,p,p,p^3))=p^{-9}(px_1^3x_2x_3+px_1^3x_2x_4+px_1^3x_3x_4-x_1^2x_2^2x_3+px_1^2x_2^2x_3
\\&\quad
-x_1^2x_2^2x_4+px_1^2x_2^2x_4+px_1^2x_2x_3^2-x_1^2x_2x_3^2-3x_1^2x_2x_3x_4+3px_1^2x_2x_3x_4
\\&\quad
+px_1^2x_2x_4^2-x_1^2x_2x_4^2-x_1^2x_3^2x_4+px_1^2x_3^2x_4+px_1^2x_3x_4^2-x_1^2x_3x_4^2
\\&\quad
+px_1x_2^3x_3+px_1x_2^3x_4-x_1x_2^2x_3^2+px_1x_2^2x_3^2+3px_1x_2^2x_3x_4-3x_1x_2^2x_3x_4
\\&\quad
+px_1x_2^2x_4^2-x_1x_2^2x_4^2+px_1x_2x_3^3-3x_1x_2x_3^2x_4+3px_1x_2x_3^2x_4+3px_1x_2x_3x_4^2
\\&\quad
-3x_1x_2x_3x_4^2+px_1x_2x_4^3+px_1x_3^3x_4+px_1x_3^2x_4^2-x_1x_3^2x_4^2+px_1x_3x_4^3
\\&\quad
+px_2^3x_3x_4+px_2^2x_3^2x_4-x_2^2x_3^2x_4-x_2^2x_3x_4^2+px_2^2x_3x_4^2+px_2x_3^3x_4
\\&\quad
+px_2x_3^2x_4^2-x_2x_3^2x_4^2+px_2x_3x_4^3)\,.
\end{align*}
These expressions are symmetrical polynomials as expected.  The written form becomes very long for higher degree $(\delta)$.  In order to be able to present intermediate results preserving the structure of these polynomials we introduce the notation of $sym_{i_1i_2i_3i_4}$ symmetrical polynomial of four variables indicating the symmetrical bundle of $x_1^{i_1}x_2^{i_2}x_3^{i_3}x_4^{i_4}$, more precisely
\begin{equation*}
sym_{i_1i_2i_3i_4} = \sum_{w\in S_4 / \,\mathrm{Stab}(x_1^{i_1}x_2^{i_2}x_3^{i_3}x_4^{i_4})} w (x_1^{i_1}x_2^{i_2}x_3^{i_3}x_4^{i_4}) \,,
\end{equation*}
where the sum is normalized by $\mathrm{Stab}(x_1^{i_1}x_2^{i_2}x_3^{i_3}x_4^{i_4})$ so the resulting coefficient is equal to 1 and $i_1\geqslant i_2\geqslant i_3\geqslant i_4\geqslant 0$.  For example,
\begin{align*}
&sym_{1000}=x_1+x_2+x_3+x_4\,,\\
&sym_{1100}=x_1x_2+x_1x_3+x_1x_4+x_2x_3+x_2x_4+x_3x_4\,,\\
&sym_{2110}=x_1^2x_2x_3+x_1^2x_2x_4+x_1^2x_3x_4+x_1x_2^2x_3+x_1x_2^2x_4+x_1x_2x_3^2\\
&\qquad+x_1x_2x_4^2+x_1x_3^2x_4+x_1x_3x_4^2+x_2^2x_3x_4+x_2x_3^2x_4+x_2x_3x_4^2\,,\\
&sym_{3333}=x_1^3x_2^3x_3^3x_4^3\,.
\end{align*}
Using this notation the previously listed examples of $\omega(\cdot)$ images become the short expressions:
\begin{align*}
&\omega(t(1,1,1,1))=1=sym_{0000}\,,\\
&\omega(t(1,1,1,p))=p^{-1}sym_{1000}\,,\\
&\omega(t(1,1,p,p))=p^{-3}sym_{1100}\,,\\
&\omega(t(1,p,p,p))=p^{-3}sym_{1110}\,,\\
&\omega(t(1,p,p,p^3))=p^{-9}(p\,sym_{3110}+(p-1)sym_{2210}+3(p-1)\,sym_{2111})\,.
\end{align*}

Finally, using formulas (\ref{eq:Ta}), (\ref{eq:T}) and (\ref{eq:omega(t)}) we obtain the images of basis Hecke operators for the symplectic group (\ref{eq:T-basis}), which we present here in $sym_{i_1i_2i_3i_4}$ notation:
\begin{equation}
\label{eq:Omega(T)}
\begin{split}
&\Omega(\Ta)=x_0\,(sym_{1111}+sym_{1110}+sym_{1100}+sym_{1000}+1)\,,\\
&\Omega(\Tb)=x_0^2 p^{-8}(
 (p-1)^2(p+1)(4p^4+3p^3+3p^2+p+1)\,sym_{1111}\\&\qquad
+p^4(p-1)(3p^2+2p+1)(sym_{2111}+sym_{1110})\\&\qquad
+p^5(p-1)(p+1)(sym_{2211}+sym_{2110}+sym_{1100})\\&\qquad
+p^7(sym_{2221}+sym_{2210}+sym_{2100}+sym_{1000}))\,,\\
&\Omega(\Tc)=x_0^2 p^{-8}(
 (p-1)(4p^4+3p^3+3p^2+p+1)\,sym_{1111}\\&\qquad
+p^2(p-1)(p^2+p+1)(sym_{2111}+sym_{1110})\\&\qquad
+p^5(sym_{2211}+sym_{2110}+sym_{1100}))\,,\\
&\Omega(\Td)=x_0^2 p^{-10}(
 (p-1)(p+1)(p^2+1)\,sym_{1111}\\&\qquad
+p^4(sym_{2111}+sym_{1110}))\,,\\
&\Omega(\Te)=x_0^2 p^{-10}\,sym_{1111}\,.
\end{split}
\end{equation}

\section{The Proof of the Formula}
\label{sec:main}

Let us return to the main subject of this article and demonstrate how the explicit result (Theorem \ref{th:main}) for generating power series was obtained.  First we find the image under spherical map of these series.  Then using formulas (\ref{eq:Omega(T)}) we can construct an equation with undetermined coefficients in order to switch from the spherical image to the expression in terms of the basis Hecke operators.

\subsection{Spherical image of Hecke power series}

Following the manipulations with the series in the book \cite{AnZh95} on page 150 we introduce the following substitutions:
\begin{align*}
\delta_2&=\delta_1+\delta_2^\prime\,,\\
\delta_3&=\delta_1+\delta_3^\prime\,,\\
\delta_4&=\delta_1+\delta_4^\prime\,,\\
\delta&=\delta_1+\delta_4^\prime+\beta\,,
\end{align*}
where $0\leqslant\delta_2^\prime\leqslant\delta_3^\prime\leqslant\delta_4^\prime\leqslant\delta^\prime$ and $\beta\geqslant 0$.  Continuing the formula (\ref{eq:Omega(D)}) using the above substitutions we obtain for $n=4$
\begin{align*}
&\Omega(\D_p(X))=\sum_{\delta=0}^{\infty}\Omega(\T(p^\delta))X^\delta =
\\&
=\sum_{\delta=0}^{\infty}\quad
 \sum_{0\leqslant\delta_1\leqslant\delta_2\leqslant\delta_3\leqslant\delta_4\leqslant\delta}\,
 p^{4\delta_1+3\delta_2+2\delta_3+\delta_4}\,
 \omega(t(p^{\delta_1},p^{\delta_2},p^{\delta_3},p^{\delta_4}))\,(x_0X)^\delta
\\&
=\sum_{\substack{\delta_1\geqslant 0, \,\beta\geqslant 0\\
       0\leqslant\delta_2^\prime\leqslant\delta_3^\prime\leqslant\delta_4^\prime}}
 (x_0X)^{\delta_1+\beta+\delta_4^\prime}
 p^{10\delta_1+3\delta_2^\prime+2\delta_3^\prime+\delta_4^\prime}
 \left(\frac{x_1x_2x_3x_4}{p^{10}}\right)^{\delta_1}
 \omega(t(1,p^{\delta_2^\prime},p^{\delta_3^\prime},p^{\delta_4^\prime}))
\\&
=\sum_{\substack{\delta_1\geqslant 0, \,\beta\geqslant 0\\
       0\leqslant\delta_2^\prime\leqslant\delta_3^\prime\leqslant\delta_4^\prime}}
 (x_0x_1x_2x_3x_4X)^{\delta_1} (x_0X)^{\beta}
 \omega(t(1,p^{\delta_2^\prime},p^{\delta_3^\prime},p^{\delta_4^\prime}))
 p^{3\delta_2^\prime+2\delta_3^\prime+\delta_4^\prime}
 (x_0X)^{\delta_4^\prime}\,.
\end{align*}
In the last formula we can separate and perform an independent summation on $\delta_1$ and $\beta$ variables.  These two series result in
\begin{equation*}
\sum_{\delta_1\geqslant 0}(x_0x_1x_2x_3x_4X)^{\delta_1}= \frac{1}{1-x_0x_1x_2x_3x_4X}
\end{equation*}
and
\begin{equation*}
\sum_{\beta\geqslant 0}(x_0X)^{\beta}= \frac{1}{1-x_0X}\,.
\end{equation*}
In the rational representation of the series $\D_p(X)=\E(X)/\F(X)$ the degree of the numerator $\E(X)$ for $n=4$ is equal to 14.  Moreover, the spherical image of the denominator $\F(X)$ is explicitly known
\begin{align*}
\Omega&(\F(X))=(1-x_0X)(1-x_0x_1X)(1-x_0x_2X)(1-x_0x_3X)(1-x_0x_4X)\\
&\times(1-x_0x_1x_2X)(1-x_0x_1x_3X)(1-x_0x_1x_4X)(1-x_0x_2x_3X)\\
&\times(1-x_0x_2x_4X)(1-x_0x_3x_4X)(1-x_0x_1x_2x_3X)(1-x_0x_1x_2x_4X)\\
&\times(1-x_0x_1x_3x_4X)(1-x_0x_2x_3x_4X)(1-x_0x_1x_2x_3x_4X)\,.
\end{align*}
Therefore we obtain
\begin{align*}
\Omega&(\E(X))=\left(\sum_{0\leqslant\delta_2^\prime\leqslant\delta_3^\prime\leqslant\delta_4^\prime}
 \omega(t(1,p^{\delta_2^\prime},p^{\delta_3^\prime},p^{\delta_4^\prime}))
 p^{3\delta_2^\prime+2\delta_3^\prime+\delta_4^\prime}
 (x_0X)^{\delta_4^\prime}\right)(1-x_0x_1X)\\
&\times(1-x_0x_2X)(1-x_0x_3X)(1-x_0x_4X)(1-x_0x_1x_2X)(1-x_0x_1x_3X)\\
&\times(1-x_0x_1x_4X)(1-x_0x_2x_3X)(1-x_0x_2x_4X)(1-x_0x_3x_4X)\\
&\times(1-x_0x_1x_2x_3X)(1-x_0x_1x_2x_4X)(1-x_0x_1x_3x_4X)(1-x_0x_2x_3x_4X)\,.
\end{align*}
In order to obtain an explicit expression for the image of the numerator $\E(X)$ we compute all
$\left(\omega(t(1,p^{\delta_2^\prime},p^{\delta_3^\prime},p^{\delta_4^\prime}))
p^{3\delta_2^\prime+2\delta_3^\prime+\delta_4^\prime}(x_0X)^{\delta_4^\prime}\right)$
up to $\delta_4^{\prime}\leqslant 14$, add them together and multiply considering only resulting powers of $X$ up to 14.  These expressions are very long, it took hours and days of processor time to compute all sums and products.  Intermediate results fill hundreds of pages of paper. However, the final result is quite short (in $sym$ notation) and it will appear separately in \cite{VaSp4}, showing some interesting properties of this polynomial (e.g. a functional equation).

\subsection{Inverting the spherical image}

In order to obtain the result of the theorem \ref{th:main} we applied the method of undetermined coefficients to each coefficient of $\Omega(\E(X))$ and $\Omega(\F(X))$.  Let us take as a reference the variable $x_0$.  In expressions for $\Omega(\E(X))$ and $\Omega(\F(X))$ this variable has the same degree as $X$ for each summand.  The expression for $\Omega(\Ta)$ (see (\ref{eq:Omega(T)})) includes the variable $x_0$ in degree 1, other images of basis Hecke operators $\Omega(\T_i(p))$ include $x_0$ in degree 2.  Therefore, to reconstruct the particularly given coefficient of degree $k$ of the polynomial $\E(X)$ or $\F(X)$ we need to construct all possible products of $\Ta$, $\Tb$, $\Tc$, $\Td$ and $\Ta$ so the resulting degree of $x_0$ in the spherical image will be equal to $k$.  For example, consider the coefficient of the degree 3 in polynomial $\E(X)$.  We computed before its image (see \cite{VaSp4}) 
\begin{align*}
&\Omega(\KE_3)= x_0^3 p^{-3}(p+1) \big( 
  p\,(sym_{3222}+sym_{3221}+sym_{3211} \\&\quad
  +sym_{3111}+sym_{2220}+sym_{2210}+sym_{2110}+sym_{1110})\\&\quad
  +(p^2+4p+1)(sym_{2222}+sym_{2221}+sym_{2211}+sym_{2111}+sym_{1111}) \big)\,.
\end{align*}
All possible products of generators of the Hecke ring having the degree 3 of $x_0$ under the spherical mapping are: $\Ta\Tb$, $\Ta\Tc$, $\Ta\Td$ and $\Te\Ta$.  Then
\begin{align*}
\Omega(\KE_3)&= K_1\Omega(\Ta)\Omega(\Tb)+K_2\Omega(\Ta)\Omega(\Tc)\\
&+K_3\Omega(\Ta)\Omega(\Td)+K_4\Omega(\Ta)\Omega(\Te)\,.
\end{align*}
Expanding these products we can construct a linear system of $K_j$ variables by comparing the coefficients of appropriate $sym$ symbols (or $x_ix_j\dotsb$ monomials).  This system resolves uniquely due to the fact that the spherical mapping constructed on basis Hecke operators is an isomorphism.  In the example above we find that $K_1=0$, $K_2=0$, $K_3=p^4(p+1)$ and $K_4=p^4(p+1)(p^2+1)(p^3-p^2+1)$.  In practice for higher degree there exist many choices of products of generators and the expansion of them becomes a not trivial task even for a computer.  For example, it took almost 80 hours of computer time to construct and resolve the linear system for the coefficient of degree 8 of the denominator.  Fortunately, there is a functional equation for coefficients of the denominator $\F(X)$ due to the symmetric structure of the spherical image polynomial:
\begin{equation*}
\KF_i=\KF_{16-i}\cdot (p^{10}\Te)^{i-8}\,,\quad i=0,\dotsc,16\,.
\end{equation*}
Therefore we used the approach of undetermined coefficient for only lower degree $\KF_i$, where $i=0,\dotsc,8$.  The same computational problem exists for the higher degree coefficient of the numerator.  To overcome the unnecessary manipulations and blind guessing of the $\T$--product combination we noticed that it is possible to lower the degree of the equation for the particular coefficient $\KE_i$ for $i>7$ by dividing this equation on factoring $\Omega(\Te)=x_0^2p^{-10}x_1x_2x_3x_4$ in appropriate degree and using the same products (with non zero coefficients) of $\Ta^{i_1}\Ta^{i_2}\Ta^{i_3}\Ta^{i_4}\Ta^{i_5}$ as for already computed $\KE_{14-i}$.

\section{Remarks}
\label{sec:rema}

\begin{rema}
The result of the Theorem \ref{th:main} is perfectly compatible with the result of the earlier work \cite{PaVaSp3}, where the same method was applied for the case of genus $g=3$.  Considering the projection from genus $g=4$ to $g=3$ corresponding to Siegel operator acting from $\Sp_4$ to $\Sp_3$ in Hecke algebra by taking $\Te_4$ to zero, $\T^{(4)}(p)$ to $\T^{(3)}(p)$, and $\T_i^{(4)}(p^2)$ to $\T_i^{(3)}(p^2)$ for $i=1, 2, 3$, we obtain the exact formula of generating power series (\ref{eq:D(n=3)}).  All formulas (\ref{eq:T-basis}) for the images of basis Hecke operators transform to the exact formulas for lower genus as well.  The spherical image $\Omega(\D_p^{(4)})$ under a projection $x_4=0$ transforms into $\Omega(\D_p^{(3)})$.  This genus lowering procedure is valid for $g=2$ and $g=1$ as well.
\end{rema}

We noticed a very interesting symmetry property within the coefficients of the spherical image of the numerator.  Knowing this relation in advance would let us to limit computation of coefficients almost in half just up to degree 7, reproducing the most time consuming higher degree coefficients using this property.

\begin{prop}
Polynomial $\Omega(\E(X))$ has the following functional relation between its coefficients $\Omega(\KE_k)$, $k=0,\dotsc ,14$:
\begin{align*}
\Omega(\KE_{14-k})&(p,x_0,x_1,x_2,x_3,x_4) = \\
&-p^{-6}(x_0^2x_1x_2x_3x_4)^{7-k}\,\Omega(\KE_{k})
 \left(\frac{1}{p},x_0x_1x_2x_3x_4,\frac{1}{x_1},\frac{1}{x_2},\frac{1}{x_3},\frac{1}{x_4}\right).
\end{align*}
\end{prop}

\begin{rema}
It is suggested that this functional relation is true for all $n$ in the following form:
\begin{equation*}
\Omega(\E)(x_0,\dotsc,x_n,X) = (-1)^{n-1}
\frac{(x_0^2x_1\dotsc x_nX^2)^{2^{n-1}-1}}{p^{\frac{n(n-1)}{2}}}
\Omega(\E)\left(\frac{1}{x_0},\dotsc,\frac{1}{x_n},\frac{p}{X}\right).
\end{equation*}
\end{rema}

For a special case of a choice Satake parameters $x_i$ the spherical image of the numerator $\E(X)$ can be considerebly simplified.
\begin{prop}
Consider the degree homomorphsm $\nu$ corresponding Satake parameters $(x_0,x_1,x_2,x_3,x_4)=(1,p,p^2,p^3,p^4)$.  Then the polynomial $\Omega(\E)$ takes the form
\begin{align*}
\Omega_\nu&(\E(X))=(1-pX)(1-p^2X)(1-p^3X)^2(1-p^4X)
\\&\times
(1+p^5X)(1-p^5X)^2(1-p^6X)^2(1-p^7X)(1-p^8X)
\\&\times
(1+pX+p^2X+2p^3X+p^4X+p^5X+2p^6X+p^7X+p^8X+p^9X^2)
\end{align*}
(compare to the similar result in genus 3 \cite{PaVaSp3}).
\end{prop}

\bigbreak
The author is very grateful to his academic advisor professor A. A. Panchishkin for posing the problem and active discussions.

\bibliographystyle{amsplain}

\end{document}